\documentclass[11pt]{amsart} 

\usepackage{amsthm,amsmath,graphicx}
\usepackage{mathrsfs} 
\usepackage[margin=3cm]{geometry}
\usepackage{hyperref} 

\usepackage[all]{xy}

\newtheorem{theorem}{Theorem}[section]
\newtheorem{conjecture}[theorem]{Conjecture} 

\theoremstyle{definition}
\newtheorem{example}[theorem]{Example}
\newtheorem{problem}[theorem]{Problem} 
\newtheorem{remark}[theorem]{Remark}

\newcommand{\Z}{\mathbb{Z}}

\newcommand{\R}{\mathbb{R}} 
\newcommand{\C}{\mathbb{C}} 
\newcommand{\PP}{\mathbb{P}}

\newcommand{\hGamma}{\widehat{\Gamma}} 
\newcommand{\hnabla}{\widehat{\nabla}} 
\newcommand{\hA}{\widehat{A}}
\newcommand{\hPhi}{\widehat{\Phi}} 
\newcommand{\hF}{\widehat{F}} 
\newcommand{\hX}{\widehat{X}} 

\newcommand{\tX}{\widetilde{X}} 

\newcommand{\cN}{\mathcal{N}} 

\newcommand{\cO}{\mathcal{O}} 
\newcommand{\cS}{\mathcal{S}} 
\newcommand{\cA}{\mathcal{A}} 
\newcommand{\cR}{\mathcal{R}} 
\newcommand{\cE}{\mathcal{E}} 

\newcommand{\cV}{\mathcal{V}} 

\newcommand{\sfC}{\mathsf{C}}

\newcommand{\frs}{\mathfrak{s}} 
\newcommand{\frL}{\mathfrak{L}} 

\newcommand{\scrD}{\mathscr{D}}

\newcommand{\ch}{\operatorname{ch}} 
\newcommand{\End}{\operatorname{End}} 
\newcommand{\id}{\operatorname{id}} 
\newcommand{\Ext}{\operatorname{Ext}} 

\newcommand{\ev}{\operatorname{ev}} 
\newcommand{\pt}{\operatorname{pt}} 
\newcommand{\QH}{\operatorname{QH}} 
\newcommand{\QDM}{\operatorname{QDM}} 
\newcommand{\ovQDM}{\overline{\QDM}}

\newcommand{\Hom}{\operatorname{Hom}}

\newcommand{\iu}{\mathtt{i}}

\newcommand{\corr}[1]{\left\langle #1 \right\rangle} 
\newcommand{\pair}[2]{\langle #1,#2 \rangle} 
\newcommand{\parfrac}[2]{\frac{\partial #1}{\partial #2}}

\numberwithin{equation}{section}

\begin{document}

\title{Gamma classes and quantum cohomology}

\author{Hiroshi Iritani}

\address{Department of Mathematics, Kitashirakawa-Oiwake-cho, Sakyo-ku, Kyoto, Japan} \email{iritani@math.kyoto-u.ac.jp}

\begin{abstract}
The $\widehat{\Gamma}$-class is a characteristic class for complex manifolds with transcendental coefficients. It defines an integral structure of quantum cohomology, or more precisely, an integral lattice in the space of flat sections of the quantum connection. We present several conjectures (the $\widehat{\Gamma}$-conjectures) about this structure, particularly focusing on the Riemann-Hilbert problem it poses. We also discuss a conjectural functoriality of quantum cohomology under birational transformations.  
\end{abstract}

\maketitle

\section{Gamma-integral structure in quantum cohomology}
We briefly review the definition of the $\hGamma$-integral structure in quantum cohomology introduced in \cite{Iritani:Int}. The corresponding rational structure was introduced independently by Katzarkov, Kontsevich and Pantev \cite{KKP:Hodge} in the framework of nc-Hodge structure. 
\subsection{Gamma class} 
\label{subsec:Gamma_class} 
Let $X$ be an almost complex manifold and let $\delta_1,\dots,\delta_n$ (with $n=\dim_\C X$) be the Chern roots of the tangent bundle, so that $c(TX) = (1+\delta_1)\cdots (1+\delta_n)$. The \emph{$\hGamma$-class} $\hGamma_X\in H^*(X;\R)$ \cite{Libgober, Lu, Iritani:Int, KKP:Hodge} is the characteristic class defined by 
\[
\hGamma_X = \Gamma(1+\delta_1) \cdots \Gamma(1+\delta_n) 
\]
where $\Gamma(1+x) = \int_0^\infty e^{-t} t^x dt$ is Euler's $\Gamma$-function. The right-hand side is expanded in symmetric power series in $\delta_1,\dots,\delta_n$ and then expressed in terms of the Chern characters $\ch_k(TX)$ as follows: 
\[
\hGamma_X = \exp\left( - \gamma c_1(X) + \sum_{k=2}^{\infty} 
(-1)^k \zeta(k) (k-1)! \ch_k(TX)\right) 
\]
where $\zeta(s) = \sum_{n=1}^\infty n^{-s}$ is the Riemann zeta function 
and $\gamma = \lim_{n\to \infty} (1+\frac{1}{2}+\cdots+\frac{1}{n}-\log n)$ 
is the Euler constant. This is a characteristic class with transcendental coefficients\footnote{It is however an algebraic (Hodge) class when $X$ is a smooth projective variety.}. 
The identity $\Gamma(1-x) \Gamma(1+x) = \pi x/\sin (\pi x)$ shows that the $\hGamma$-class can be thought of as a ``square root'' of the $\hA$-class, i.e.
\begin{equation} 
\label{eq:hGamma_hA}
\hGamma_X \cdot \hGamma_X^* = (2\pi \iu)^{\deg/2} \hA_X 
\end{equation}
where $\hGamma_X^* := (-1)^{\deg/2} \hGamma_X$ denotes the dual $\hGamma$-class. We note that $\hA_X$ depends only on the underlying topological manifold whereas $\hGamma_X$ depends on an almost complex structure on it. 
The identity \eqref{eq:hGamma_hA} suggests a relationship between the $\hGamma$-class and the Atiyah-Singer index theorem. In fact, we can interpret $\hGamma_X$ as (a regularization of) the inverse $S^1$-equivariant Euler class of the positive normal bundle $\cN_+$ of the set $X$ of constant loops in the free loop space $LX$ (see \cite{Lu} \cite[Appendix A]{GGI}), i.e.~ 
\begin{equation} 
\label{eq:Euler_positive_normal}
\frac{1}{e_{S^1}(\cN_+)} = \frac{1}{\prod_i \prod_{k>0} (\delta_i + kz)} \sim (2 \pi)^{-n/2} z^{(n-\deg)/2} z^{c_1(X)} \hGamma_X
\end{equation}
where $z$ is a generator of the $S^1$-equivariant cohomology of a point. This is reminiscent of the loop space heuristics of the index theorem by Atiyah and Witten \cite{Atiyah:circular}, where the $\hA$-class is interpreted as the inverse Euler class $e_{S^1}(\cN)^{-1}$ of the normal bundle $\cN$ itself. 
\[
\frac{1}{e_{S^1}(\cN)} =\frac{1}{e_{S^1}(\cN_-) e_{S^1}(\cN_+)} 
= \frac{1}{\prod_i \prod_{k\neq 0} (\delta_i + kz)} \sim  \left(\frac{z}{2\pi\iu}\right)^{n-(\deg/2)} \hA_X.
\] 
Since $\cN_+$ corresponds to infinitesimal (pseudo-)holomorphic loops, the $\hGamma$-class can be thought of as the \emph{localization contribution from constant loops} in symplectic Floer theory. 

\subsection{Quantum cohomology D-modules}
\label{subsec:QDM} 
Let $X$ be a smooth projective variety (or a compact symplectic manifold) 
and let $H^*(X)$ denote the cohomology group with complex coefficients. 
The quantum cohomology $\QH^*(X) = (H^*(X),\star_\tau)$ of $X$ is a family of supercommutative product structures $\star_\tau$ on $H^*(X)$ parametrized by $\tau \in H^*(X)$. The quantum product $\star_\tau$ is defined by 
\[
(\alpha\star_\tau \beta,\gamma) = \sum_{d\in H_2(X,\Z), k\ge 0} 
\corr{\alpha,\beta,\gamma,\tau,\dots,\tau}_{0,k+3,d} \frac{Q^d}{k!} 
\]
for $\alpha,\beta,\gamma \in H^*(X)$. Here $(\alpha,\beta) = \int_X \alpha \cup \beta$ is the Poincar\'e pairing and $\corr{\alpha_1,\dots,\alpha_k}_{0,k,d}$ denotes the genus-zero, $k$-point, degree $d$ Gromov-Witten invariants. Strictly speaking, we should treat the odd degree part of $\tau$ as anti-commuting variables and view the parameter space $H^*(X)$ as a supermanifold. For the most part of this paper, we shall restrict the parameter $\tau$ and elements of quantum cohomology to the even part of the cohomology group and write $H^*(X)$ for the even part (see Remark \ref{rem:odd} for the odd part). 

In the above formula, we introduced the Novikov variable $Q$ to ensure the adic convergence of $\star_\tau$. The divisor equation shows that, if we decompose $\tau = \sigma + \tau'$ with $\sigma \in H^2(X)$ and $\tau' \in H^{\neq 2}(X)$, 
\[
(\alpha\star_\tau \beta,\gamma) = \sum_{d\in H_2(X,\Z), k\ge 0} 
\corr{\alpha,\beta,\gamma,\tau',\dots,\tau'}_{0,k+3,d} \frac{e^{\pair{\sigma}{d}} Q^d}{k!}. 
\]
Thus the quantum product can be expanded in a power series in $\tau'$ and $e^{\sigma}$ and approaches the cup product in the following \emph{large-radius limit}: 
\begin{equation} 
\label{eq:LRL}
\tau' \to 0, \quad 
\text{$e^{\pair{\sigma}{d}} \to 0$ for all effective classes $d\neq 0$}.  
\end{equation} 
Hereafter we shall always specialize the Novikov variable $Q$ to $1$ and assume that $\star_\tau|_{Q=1}$ (which we shall write as $\star_\tau$)  is convergent in a neighbourhood $U$ of the large radius limit. 

The quantum cohomology defines the structure of a Frobenius manifold \cite{Dubrovin:ICM} on the convergence domain $U\subset H^*(X)$. Specifically, it defines a meromorphic flat connection $\nabla$ on the trivial bundle $F=H^*(X)\times (U \times \C) \to (U \times \C)$, called the \emph{quantum connection} or the \emph{Dubrovin connection}. It is defined by the formulae 
\begin{align*} 
\nabla_{\partial/\partial \tau^i} &= \parfrac{}{\tau^i} + \frac{1}{z} (\phi_i\star_\tau) \\
\nabla_{z \partial/\partial z} &= z \parfrac{}{z} - \frac{1}{z} (E\star_\tau) + \mu 
\end{align*} 
where $(\tau,z) \in U \times \C$ denotes a point on the base and $\{\tau^i\}$ are linear coordinates dual to a homogeneous basis $\{\phi_i\}$ of $H^*(X)$ so that $\tau = \sum_i \tau^i \phi_i$. The section $E\in \cO(F)$ is the \emph{Euler vector field} given by 
\[
E = c_1(X) + \sum_{i} \left(1 -\frac{\deg \phi_i}{2}\right) \tau^i \phi_i 
\]
and $\mu\in \End(H^*(X))$ is the \emph{grading operator} defined by $\mu (\phi_i) = (\frac{\deg \phi_i}{2} -\frac{n}{2}) \phi_i$. The connection $\nabla$ has poles of order two along $z=0$ and is possibly \emph{irregular singular} there. On the other hand, it has logarithmic poles (and is therefore \emph{regular singular}) along $z=\infty$. The connection $\nabla$ is compatible with the Poincar\'e pairing in the sense that the $z$-sesquilinear pairing 
\begin{equation} 
\label{eq:sesqui_pairing} 
(-1)^*\cO(F) \otimes \cO(F) \to \cO, \qquad s_1 \otimes s_2 \mapsto (s_1,s_2) = \int_X s_1 \cup s_2  
\end{equation} 
is flat for $\nabla$, where $(-1) \colon U\times \C \to U\times \C$ is the map sending $(\tau,z)$ to $(\tau,-z)$. 

The \emph{quantum (cohomology) D-module} is the tuple $\QDM(X) = (\cO(F), \nabla, (\cdot,\cdot))$ consisting of the locally free sheaf $\cO(F)$ over $U\times\C$, the quantum connection $\nabla$ and the pairing in \eqref{eq:sesqui_pairing}. 
The D-module approach to quantum cohomology has been proposed by Givental \cite{Givental:homological} and Guest \cite{Guest:D-mod}. 

\subsection{Gamma-integral structure} 
The $\hGamma$-integral structure is an integral lattice in the space of flat sections for the quantum connection $\nabla$. We have a fundamental solution for $\nabla$-flat sections of the form $L(\tau,z) z^{-\mu} z^{c_1(X)}$ with an $\End(H^*(X))$-valued function $L(\tau,z)$ uniquely characterized by the following asymptotic condition at the large radius limit \eqref{eq:LRL}: 
\[
L(\tau,z) = (\id + O(e^{\sigma},\tau')) e^{-\sigma/z}.
\]
Here $z^{-\mu} z^{c_1(X)} = \exp(-\mu \log z) \exp(c_1(X) \log z)$ is an $\End(H^*(X))$-valued function on the universal cover of $\C^\times$. 
Given a basis $\{\phi_i\}$ of $H^*(X)$, $\{L(\tau,z) z^{-\mu} z^{c_1(X)} \phi_i\}$ gives a basis of $\nabla$-flat sections. 
Explicitly, $L(\tau,z)$ is given in terms of gravitational descendants as follows (see \cite{Dubrovin:2D}, \cite{Givental:equivariant}): 
\begin{equation} 
\label{eq:fundsol}
L(\tau,z) \phi_i = e^{-\sigma/z}\phi_i - \sum_j \sum_{(d,k) \neq (0,0)} 
\corr{\frac{e^{-\sigma/z} \phi_i}{z+\psi}, \tau',\dots,\tau', \phi^j}_{0,k+2,d} 
\frac{\phi_j}{k!} e^{\pair{\sigma}{d}} 
\end{equation} 
where $\psi$ denotes the universal cotangent class at the first marking, $1/(z+\psi)$ should be expanded as $\sum_{k\ge 0} z^{-k-1} (-\psi)^k$, and $\{\phi_j\}$, $\{\phi^j\}$ are mutually dual bases of $H^*(X)$ such that $\int_X \phi_i \cup \phi^j = \delta_i^j$. 

Let $\cS(X)$ denote the $\C$-vector space of multi-valued flat sections of $(F,\nabla)$ over $U\times \C^\times$, i.e.~flat sections over the universal cover $U\times \widetilde{\C^\times}$. Let $K(X) = K^0_{\textrm{top}}(X)$ denote the $K$-group of topological complex vector bundles and define a map $\frs \colon K(X) \to \cS(X)$ as 
\begin{equation} 
\label{eq:framing}
\frs(V)(\tau,z) = L(\tau,z) z^{-\mu} z^{c_1(X)} 
\left( (2\pi)^{-n/2} \hGamma_X (2\pi\iu)^{\deg/2} \ch(V)\right). 
\end{equation} 
The factor $(2\pi)^{-n/2} z^{-\mu} z^{c_1(X)}\hGamma_X$ also appears in \eqref{eq:Euler_positive_normal} as a regularization of $e_{S^1}(\cN_+)^{-1}$.  The \emph{$\hGamma$-integral structure} is the integral lattice of $\cS(X)$ given as the image of the map $\frs$.  

By the compatibility between $\nabla$ and the Poincar\'e pairing, we have a non-degenerate (not necessarily symmetric or anti-symmetric) pairing $[\cdot,\cdot) \colon \cS(X)\otimes \cS(X) \to \C$ defined by 
\begin{equation} 
\label{eq:[)}
[s_1,s_2) = (s_1(\tau,e^{-\pi \iu}z), s_2(\tau,z)) 
\end{equation} 
for $s_1,s_2 \in \cS(X)$. The property \eqref{eq:hGamma_hA} of the $\hGamma$-class and the Atiyah-Singer index theorem (or Hirzebruch-Riemann-Roch theorem) show that $\frs$ respects the pairing 
\[
[\frs(V_1), \frs(V_2)) = \chi(V_1,V_2)
\]
where $\chi(V_1,V_2)\in \Z$ is the $K$-theoretic push-forward of $V_1^\vee\otimes V_2$ to a point (the index of a Dirac operator; it is $\sum_{i\ge 0} (-1)^i \dim \Ext^i(V_1,V_2)$ if $V_1$, $V_2$ are holomorphic vector bundles). 
The $\hGamma$-integral structure is monodromy-invariant in the sense that 
\begin{align*} 
\frs(V)(\tau - 2\pi \iu c_1(L), z) &= \frs(V\otimes L)(\tau,z) \\ 
\frs(V)(\tau, e^{-2\pi \iu}z)& = \frs(V \otimes \omega_X[n]) (\tau,z)
\end{align*} 
where $L$ is a (topological) line bundle on $X$ and $\omega_X[n]=(-1)^n \omega_X$ is the canonical line bundle $\omega_X$ shifted by $n$, corresponding to the Serre functor of the derived category. 

\begin{remark}
The $\hGamma$-integral structure can be defined more generally for orbifolds \cite{Iritani:Int}. 
\end{remark} 

\begin{remark} 
\label{rem:odd} 
We can generalize the $\hGamma$-integral structure including the odd part of the quantum cohomology, using $K^*(X) = K^0_{\textrm{top}}(X) \oplus K^1_{\textrm{top}}(X)$ instead of $K^0_{\textrm{top}}(X)$. We still restrict the parameter $\tau$ to lie in the even part, but consider flat sections taking values in the full cohomology group. The formula \eqref{eq:framing} makes sense for all $V \in K^*(X)$ when we choose a square root $\sqrt{2 \pi \iu}$ and use the Chern character $\ch\colon K^*(X) \to H^*(X)$ of Atiyah-Hirzebruch \cite{Atiyah-Hirzebruch}. The resulting map $\frs_X \colon K^*(X)/\textrm{tors} \to \cS(X)$ then has the advantage that it is natural with respect to the Cartesian product, i.e.~$\frs_{X\times Y} = \frs_X \otimes \frs_Y$ (under the K\"unneth isomorphism). Interestingly, we have 
\begin{align*} 
[\frs(\alpha_1), \frs(\alpha_2)) & = -\iu \chi(\alpha_1,\alpha_2) \in \iu \Z \qquad && \text{for $\alpha_1,\alpha_2 \in K^1(X)$}\\
\frs(\alpha)(\tau,e^{-2\pi\iu} z) & = (-1)^{\deg \alpha} \frs(\alpha\otimes \omega_X[n])(\tau,z) \qquad && \text{for $\alpha \in K^*(X)$}  
\end{align*} 
where $\chi(\alpha_1,\alpha_2)\in \Z$ is the $K$-theoretic push-forward of $\alpha_1^\vee \cdot \alpha_2$ to a point as before; the dual element $\alpha_1^\vee$ here is defined via the isomorphism $K^1(X) \cong K^{-1}(X) \cong \tilde{K}^0(S^1\wedge X^+)$ (see \cite{Atiyah-Hirzebruch}) and the usual duality in $K^0$. 
\end{remark} 

\begin{remark}[mirror symmetry] 
The $\hGamma$-integral structure had (implicitly) appeared for a long time in the study of mirror symmetry before it was defined in \cite{Iritani:Int, KKP:Hodge}. Under mirror symmetry of Calabi-Yau manifolds, the quantum differential equation corresponds to the Picard-Fuchs differential equations satisfied by periods of the mirror family, and we can partially see the $\hGamma$-class in the asymptotics of periods near the large-complex structure limit. Libgober \cite{Libgober} introduced the (inverse) $\hGamma$-class based on the observation of Hosono et.~al.~\cite{HKTY} that certain combinations of Chern numbers and $\zeta$ values appear in solutions of the mirror Picard-Fuchs equations. Hosono \cite{Hosono:central} stated a conjecture equating periods of mirrors of complete intersections with explicit hypergeometric series and the $\hGamma$-class is hidden in the series. We also refer the reader to  \cite{Horja, Borisov-Horja:FM,vEvS} for related works. It has been checked in a number of cases that the $\hGamma$-integral structure corresponds to a natural integral structure on the mirror side \cite{Iritani:Int, Iritani:periods}. Regarding the compatibility with mirror symmetry, an approach based on the SYZ picture and tropical geometry has been proposed in \cite{AGIS} recently.   
\end{remark} 

\section{Gamma conjectures}
In this section we review the $\hGamma$-conjectures I, II discussed by Galkin, Golyshev and the author \cite{GGI}, and their generalization by Sanda and Shamoto \cite{Sanda-Shamoto:analogue}. The $\hGamma$-conjectures can be understood as the compatibility between the Betti (real, rational, or integral) structure and the Stokes structure, discussed by Hertling and Sevenheck \cite{Hertling-Sevenheck:nilpotent} in the context of TERP structure and by Katzarkov, Kontevich and Pantev \cite{KKP:Hodge} in the context of nc-Hodge structure. The Gamma conjecture II also refines Dubrovin's conjecture \cite{Dubrovin:ICM}. 

\subsection{Gamma conjecture I} 
\label{subsec:Gamma_I} 
The $\hGamma$-conjecture I is specifically about quantum cohomology of Fano manifolds. It roughly speaking says that we can know the topology (the $\hGamma$-class) of a Fano manifold by counting rational curves on it. In view of \eqref{eq:hGamma_hA}, we may view it as a ``square root'' of the index theorem. 

Let $X$ be a Fano manifold and let $J_X(\tau,z)$ be the (small) $J$-function defined as 
\[
J_X(\tau,z) = e^{\tau/z} \left(1 + \sum_i \sum_{d\in H_2(X,\Z), d\neq 0} 
e^{\pair{\tau}{d}}\corr{\frac{\phi^i}{z(z-\psi)}}_{0,1,d} \phi_i  \right) 
\]
where $\tau \in H^2(X)$. This is a cohomology-valued function which is convergent for all $(\tau,z)\in H^2(X)\times \C^\times$ (this follows from the Fano assumption). We can also write this as $J_X(\tau,z) = L(\tau,z)^{-1} 1$ using the fundamental solution $L$ in \eqref{eq:fundsol}; hence $J_X(\tau,z)$ gives a solution of the quantum D-module along the $\tau$-direction. 

\begin{conjecture}[$\hGamma$-conjecture I] 
\label{conj:Gamma_I} 
For a Fano manifold $X$, we have the equality 
\[
 [ \hGamma_X ] = \lim_{t\to +\infty} [ J_X(c_1(X) \log t, 1) ]
\]
in the projective space $\PP(H^*(X))$ of cohomology. 
\end{conjecture} 
This has been proved for the projective spaces, type A Grassmannians \cite{Golyshev:deresonating,GGI} and Fano threefolds of Picard rank one \cite{Golyshev-Zagier}. The $\hGamma$-conjecture I for Fano toric manifolds or complete intersections in them follows if these spaces satisfy certain conditions related to Conjecture $\cO$ \cite{GI}. The $\hGamma$-conjecture I is also compatible with taking hyperplane sections, i.e.~if a Fano manifold $X$ satisfies the $\hGamma$-conjecture I and if $Y\subset X$ is a hypersurface in the linear system $|L|$ with $L$ proportional to $-K_X$, $Y$ satisfies the $\hGamma$-conjecture I \cite[Theorem 8.3]{GI}.   
\begin{example} 
The $J$-function of $\PP^n$ is given by 
\[
J_{\PP^n}( c_1(\PP^n) \log t,z) = \sum_{d=0}^\infty 
\frac{t^{(n+1)(d+p/z)}}{\prod_{k=1}^d (p+kz)^{n+1}}
\]
where $p$ is the hyperplane class. Setting $z=1$ and fixing $t>0$, we find that the $d$th summand 
\[
\frac{t^{(n+1)(d+p)}}{\prod_{k=1}^d (p+k)^{n+1}} 
= \left( 
\frac{t^d (t/d)^{p}}{d!} \right)^{n+1} 
\left( \frac{e^{(\log d-(1+\frac{1}{2}+ \cdots + \frac{1}{d}))p}}{\prod_{k=1}^d((1+\frac{p}{k})e^{-p/k})} 
\right)^{n+1} 
\]
has a strong peak approximately when $d$ is close to $t$. We can guess from this that the limit of $\C J_{\PP^n}(c_1(\PP^n) \log t,1)$ in the projective space should be the line generated by 
\[
\lim_{d\to \infty} 
\left( \frac{e^{(\log d-(1+\frac{1}{2}+ \cdots + \frac{1}{d}))p}}{\prod_{k=1}^d((1+\frac{p}{k})e^{-p/k})} 
\right)^{n+1} 
= \Gamma(1+p)^{n+1} = \hGamma_{\PP^n}. 
\]
\end{example} 

\begin{remark} In Givental's heuristic calculation of the $J$-function \cite{Givental:homological}, the $d$th summand $\prod_{k=1}^d (p+kz)^{-n-1}$ of $J_{\PP^n}(\tau,z)$ appears as the localization contribution from constant loops in the polynomial loop space (quasi-maps' space) of degree $d$, and hence it can be viewed as the degree-$d$ truncation of $e_{S^1}(\cN_+)^{-1}$ appearing in \S\ref{subsec:Gamma_class}. In view of the loop space interpretation of the $\hGamma$-class, this gives a geometric explanation for the $\hGamma$-conjecture I in this case. In general, the degree $d$ term of $J_X(\tau,z)$ arises from a localization contribution of an integral over the graph space $G_d = \overline{M}_{0,0}(X\times \PP^1, (d,1))$, the moduli space of genus-zero stable maps to $X\times \PP^1$ of degree $(d,1)$, equipped with the $\C^\times$-action induced by the $\C^\times$-action on $\PP^1$. Let $G_d^\circ \subset G_d$ denote the open subset consisting of stable maps which are genuine graphs near $\infty\in \PP^1$, i.e.~do not contain components contained in $X\times \{\infty\}$. Then $G_d^\circ$ is preserved by the $\C^\times$-action and has $F_d = \overline{M}_{0,1}(X,d)$ as the fixed locus. Writing $\ev_\infty \colon G_d^\circ \to X$ for the evaluation map at $\infty \in \PP^1$, we have 
\[
\int_{G_d^\circ} \ev_\infty^*\alpha = 
\int_{[F_d]_{\textrm{vir}}} \frac{\ev_1^*\alpha}{z(z-\psi)} = (\text{degree $d$ term of $J_X$}, \alpha)  
\]
where we defined the integral over the improper space $G_d^\circ$ using (virtual)  equivariant localization. In the case where $X= \PP^n$, Givental gave a birational morphism from $G_d^\circ$ to the polynomial loop space and justified his heuristic calculation (see \cite[Main Lemma]{Givental:equivariant}). Therefore the $\hGamma$-conjecture I can be viewed as the statement that \emph{$G_d^\circ$ approximates the (positive) loop space of $X$ as $d\to \infty$} in a suitable sense. 
\end{remark}


\begin{remark} 
A discrete version of the limit in Conjecture \ref{conj:Gamma_I} had been also studied (before the formulation of the $\hGamma$-conjecture) and called \emph{Ap\'ery limits}, in view of the connection to Ap\'ery's proof of the irrationality of $\zeta(2)$ and $\zeta(3)$, see Galkin \cite{Galkin:Apery} and Golyshev \cite{Golyshev:deresonating}.
\end{remark} 

\subsection{Gamma conjecture I in terms of flat sections} 

We restate $\hGamma$-conjecture I in terms of $\nabla$-flat sections over the $z$-plane, in order to explain the relationship with $\hGamma$-conjecture II in the following section. We start with Conjecture $\cO$. 
\begin{conjecture}[Conjecture $\cO$] 
Let $X$ be a Fano manifold and let $T$ denote the maximal norm of the eigenvalues of the quantum multiplication $(E\star_0) = (c_1(X) \star_0)$ at $\tau=0$. 
Then $T$ is a simple eigenvalue of $(c_1(X)\star_0)$, that is, an eigenvalue whose multiplicity in the characteristic polynomial is one.  
\end{conjecture} 
Here we omitted part (2) of Conjecture $\cO$ in \cite[Definition 3.1.1]{GGI} since we do not need it. Conjecture $\cO$ is a consequence of the Perron-Frobenius theorem if $(c_1(X)\star_0)$ is represented by an irreducible nonnegative matrix. Cheong and Li \cite{Cheong-Li} proved Conjecture $\cO$ for homogeneous spaces $G/P$ using the Perron-Frobenius theorem. 

Eigenvalues of $(c_1(X)\star_0)$ are closely related to asymptotics of flat sections for $\nabla|_{\tau=0}$ as $z\to 0$. The flatness equation reads 
\[
\left(z\parfrac{}{z} - \frac{1}{z} c_1(X)\star_0 + \mu \right) s(z) = 0. 
\]
For each eigenvector $\Psi$ of $(c_1(X)\star_0)$ with eigenvalue $u$, we expect that there should be a flat section $s(z)$ with asymptotics $\sim e^{-u/z}\Psi$ as $z\to 0$. 
Define a vector space $\cA$ as 
\[
\cA = 
\left\{s \colon \R_{>0} \to H^*(X)\, \middle| \,
\begin{array}{l} 
\text{$s(z)$ is flat for $\nabla_{z\partial_z}|_{\tau=0}$} \\  
\text{$e^{T/z} s(z)$ is at most of polynomial growth as $z\to +0$}
\end{array} 
\right\}
\]
Assuming Conjecture $\cO$ for $X$, we can prove that $\cA$ is one-dimensional and that, for any $s(z) \in \cA$, $e^{T/z}s(z)$ converges to a $T$-eigenvector of $(c_1(X)\star_0)$ as $z\to +0$. The original formulation of the $\hGamma$-conjecture I in \cite{GGI} was as follows. 
\begin{conjecture}[$\hGamma$-conjecture I: another form] 
The space $\cA$ is generated by $\frs(\cO)|_{\tau=0}$.  
\end{conjecture} 
It is equivalent to Conjecture \ref{conj:Gamma_I} in \S \ref{subsec:Gamma_I} under Conjecture $\cO$ \cite[Corollary 3.6.9]{GGI} and can be viewed as a dual formulation. We have a gauge equivalence between the connections $\nabla|_{\tau=0}$ and $\nabla|_{z=1, \tau = c_1(X) \log t}$, that is, $z^{\mu} \left(\nabla|_{\tau=0}\right) z^{-\mu} = \nabla|_{z=1, \tau=c_1(X) \log t}$ under the identification $t= z^{-1}$.Therefore flat sections over $\{\tau=0\}\times \C$ and the solution $J_X(c_1(X)\log t,1)$ are dual to each other. While the space $\cA$ consists of flat sections with most rapid decay ($\sim e^{-T/z}$), the $t\to +\infty$ limit of the $J$-function detects its most rapidly growing component. In fact, we can show that the $J$-function has the following asymptotics:  
\[
J_X(c_1(X) \log t, 1) = C t^{-n/2} e^{Tt} \left(\hGamma_X +O(t^{-1}) \right) \quad \text{as $t\to +\infty$}
\]
for some $C \in \C$, under Conjecture $\cO$ and $\hGamma$-conjecture I (see \cite[Proposition 3.8]{GI}). 

\subsection{Gamma conjecture II} 
In this section we assume that the quantum product $\star_\tau$ is semisimple\footnote
{Under the semisimlicity assumption, $X$ has no odd cohomology classes and, if moreover $X$ is a smooth projective variety, $H^*(X)$ is necessarily of Hodge-Tate type, i.e.~$H^{p,q}(X) = 0$ for $p\neq q$, see \cite{HMT}.} at some $\tau=\tau_0\in H^*(X)$, i.e.~$(H^*(X),\star_{\tau_0})$ is isomorphic to the direct sum of $\C$ as a ring. We do not need to assume that $X$ is Fano. Let $\psi_1,\dots,\psi_N\in H^*(X)$ denote an idempotent basis such that $\psi_i \star_\tau \psi_j = \delta_{ij} \psi_i$ and let $u_1,\dots,u_N\in \C$ be the eigenvalues of $(E\star_\tau)$ such that $E\star_\tau \psi_i = u_i \psi_i$; here $\psi_i$ and $u_i$ are analytic functions of $\tau$ defined in a neighbourhood of $\tau=\tau_0$. The functions $\{u_i\}$ give a local coordinate system near $\tau_0$ called the \emph{canonical coordinates} \cite{Dubrovin:2D,Dubrovin:Painleve}. We write $\Psi_i = (\psi_i,\psi_i)^{-1/2}\psi_i$ for the normalized idempotent basis, which is unique up to sign. Choose a phase $\phi \in \R$ such that $e^{\iu\phi}\notin \R_{>0}(u_{i,0} - u_{j,0})$ for all $i,j$, where $u_{i,0}$ is the value of $u_i$ at $\tau_0$; such a phase $\phi$ is said to be \emph{admissible}.
We have a basis $(y^\phi_1(\tau,z),\dots,y^\phi_N(\tau,z))$ of $\nabla$-flat sections defined in a neighbourhood of $\tau=\tau_0$ and $\arg z = \phi$ with the following property: 
\[
e^{u_i/z} y^\phi_i(\tau,z) \to \Psi_i \quad \text{as $z \to 0$ along the angular sector $|\arg z -\phi|<\pi + \epsilon$}
\]
for some $\epsilon>0$, see \cite[Proposition 2.5.1]{GGI}. 
\begin{conjecture}[$\hGamma$-conjecture II: a topological form] 
\label{conj:Gamma_II} 
Suppose that the quantum product $\star_\tau$ of $X$ is semisimple at some $\tau_0\in H^*(X)$ and let $\phi\in \R$ be an admissible phase for the eigenvalues of $(E\star_{\tau_0})$. 
There exist $K$-classes $\cE^\phi_1,\dots,\cE^\phi_N \in K(X)$ such that $y_i^\phi(\tau,z) = \frs(\cE_i^\phi)(\tau,z)$ in a neighbourhood of $\tau=\tau_0$ and $\arg z = \phi$.
\end{conjecture} 
This refines part (3) of Dubrovin's conjecture \cite[Conjecture 4.2.2]{Dubrovin:ICM} concerning the central connection matrix. It has been proved for type A Grassmannians \cite{GGI}, Fano toric manifolds \cite{Fang-Zhou:GammaII_toric} and quadric hypersurfaces \cite{Hu-Ke:GammaII_quadric}. 
The $\hGamma$-conjecture I can be viewed as a special case of the $\hGamma$-conjecture II when $\tau=0$ and $\phi=0$. 

The flat sections $y_i^\phi(\tau,z)$ depend on the choice of a phase $\phi$ (or more precisely on a chamber of admissible phases) whereas their asymptotic expansions as $z\to 0$ do not. This is the so-called Stokes phenomena. The \emph{Stokes matrix} $S=(S_{ij})$ is a transition matrix between the flat sections associated with opposite directions: it is given by $y_j^{\phi}(\tau,z) = \sum_{i=1}^N y_i^{\phi+\pi}(\tau,z) S_{ij}$ with $\arg z = \phi + \frac{\pi}{2}$. 
It can be given in terms of the bilinear form in \eqref{eq:[)} and then as the Euler matrix of $\{\cE_i^\phi\}$: 
\begin{equation*} 
S_{ij} = (y_i^\phi, y_j^\phi] = \chi(\cE_i^\phi, \cE_j^\phi). 
\end{equation*} 
This corresponds to part (2) of Dubrovin's conjecture saying that the Stokes matrix is integral and is given by the Euler pairing. If follows from the fact that the asymptotics $y_i^\phi \sim e^{-u_i/z} \Psi_i$ holds over a sector of angle $>\pi$ that the Stokes matrix is upper-triangular
\[
S_{ij} = \chi(\cE_i^\phi,\cE_j^\phi) = \begin{cases} 
1 & \text{if $i=j$;} \\
0 & \text{if $\Im (e^{-\iu\phi} u_i) \le \Im(e^{-\iu\phi} u_j)$ and $i\neq j$.}
\end{cases} 
\]
\begin{remark} In \cite{GGI}, the $\hGamma$-conjecture II was stated for Fano manifolds which have semisimple quantum cohomology and full exceptional collections in $D^b(X)$. It is moreover conjectured that $\{\cE_i^\phi\}$ should lift to a full exceptional collection. We drop these assumptions/conclusions to emphasize a  topological nature of the $\hGamma$-conjecture. 
\end{remark} 
\begin{remark} Dubrovin \cite{Dubrovin:Strasbourg} also formulated a conjecture similar to the $\hGamma$-conjecture II. See Cotti, Dubrovin and Guzzetti \cite{CDG:helix} for the formulation. 
\end{remark} 

\begin{example}
For $X=\PP^n$, the corresponding exceptional collection is $\{\cO, \cO(1),\dots,\cO(n)\}$ at some $\tau$ \cite{GGI}. The collection at $\tau=0$ is given explicitly in \cite{CDG:helix}. 
\end{example} 

\begin{remark} 
Suppose that $X$ is Fano and is mirror to a Landau-Ginzburg model $f \colon Y \to \C$. It is expected that the idempotent $\psi_i$ corresponds to a non-degenerate critical point $c_i$ of $f$ such that the corresponding eigenvalue $u_i$ equals $f(c_i)$. The critical point $c_i$ can associate a Lefschetz thimble $\frL_i^\phi$ extending in the direction of $e^{\iu\phi}$, which gives an exceptional object in the Fukaya-Seidel category of $(Y,f)$. The object in $D^b(X)$ corresponding to $\frL_i^\phi$ under homological mirror symmetry should give the class $\cE_i^\phi$. 
\end{remark}

\begin{remark} 
The $\hGamma$-conjecture II concerns the connection problem between flat sections $y_i^\phi$ characterized by the asymptotics at the irregular singular point $z=0$ and flat sections $\frs(V)$ normalized at the regular singular point $z=\infty$. The connection matrix of flat sections (with respect to a fixed basis) is called the \emph{central connection matrix} by Dubrovin \cite{Dubrovin:ICM}; in the formalism of the $\hGamma$-integral structure, it corresponds to the basis $\{\cE_i^\phi\}$ of the $K$-group. As discussed by Dubrovin, when $\tau$ and $\phi$ vary, the basis $\{y_i^\phi\}$ of flat sections (and hence $\{\cE_i^\phi\}$) changes discontinuously by the action of the braid group in $N$ strands. Suppose that we ordered flat sections $\{y_i^\phi\}$ in such a way that $\Im(e^{-\iu\phi} u_1) \ge \Im(e^{-\iu\phi} u_2) \ge \cdots \ge \Im(e^{-\iu\phi} u_N)$. The braid group action is generated by the following right mutations and their inverses (which are the actions of simple braids): 
\[
(\cE_1,\dots,\cE_i, \cE_{i+1}, \dots, \cE_N)
\mapsto (\cE_1, \dots, \cE_{i+1}, \cE_i - \chi(\cE_i,\cE_{i+1}) \cE_{i+1},\dots,\cE_N)   
\]
where we suppressed $\phi$ to simplify the notation. 
This transformation happens when the eigenvalue $u_i$ crosses behind $u_{i+1}$ towards direcion $e^{\iu\phi}$. As Dubrovin observed, this is consistent with mutations of exceptional collections in the derived category \cite{Bondal-Polishchuk}. 
\end{remark} 

\subsection{Conjecture of Sanda and Shamoto} 
\label{subsec:Gamma_III} 
Sanda and Shamoto \cite{Sanda-Shamoto:analogue} proposed a generalization of the  $\hGamma$-conjecture II to the case where quantum cohomology is not necessarily semisimple and called it \emph{Dubrovin-type conjecture}. Their formulation involves derived category of coherent sheaves and Hochschild homology, but here we give a topological formulation that has been proposed by Sergey Galkin \cite{Galkin:private}. This formulation makes sense for any compact symplectic manifolds.  

We fix a parameter $\tau \in H^*(X)$ in the convergence domain of the quantum product. We consider the restriction $\QDM(X)_\tau$ of the quantum D-module from \S\ref{subsec:QDM} to $\{\tau\} \times \C$ and write  $\ovQDM(X)_\tau=\QDM(X)_{\tau,0}\otimes_{\C\{z\}}\C[\![z]\!]$ for the restriction to the formal neighbourhood of $z=0$, where $\QDM(X)_{\tau,0}$ denotes the germ of $\QDM(X)_\tau$ at $z=0$. We say that the quantum connection at $\tau$ is of \emph{exponential type}\footnote{We follow the terminology in \cite{KKP:Hodge}; it was called ``\emph{require no ramification}'' in \cite{Hertling-Sevenheck:nilpotent}.} if we have the following formal decomposition (see \cite[Lemma 8.2]{Hertling-Sevenheck:nilpotent}): 
\begin{equation} 
\label{eq:formal_decomp} 
\hPhi \colon \ovQDM(X)_\tau 
\cong \bigoplus_{u\in \sfC} 
(e^{u/z} \otimes \cR_u) \otimes_{\C\{z\}} \C[\![z]\!] 
\end{equation} 
where we disregard the pairing on $\ovQDM(X)_\tau$ momentarily, $\sfC$ denotes the set of distinct eigenvalues of $(E\star_\tau)$, $e^{u/z}$ denotes the rank-one connection $(\C\{z\}, d + d(u/z))$ and $\cR_u$ is a free $\C\{z\}$-module equipped with a regular singular connection (whose pole order at $z=0$ is at most two). In this decomposition, each ``regular singular piece'' $\cR_u$ is unique up to isomorphism. This decomposition is automatically orthogonal with respect to the Poincar\'e pairing \eqref{eq:sesqui_pairing} and hence each piece $\cR_u$ inherits a non-degenerate $z$-sesquilinear pairing $(\cdot,\cdot)_u\colon (-1)^*\cR_u \otimes \cR_u \to \C\{z\}$. Hereafter we assume that the quantum connection is of exponential type: this assumption is natural from a mirror symmetry point of view. 

We choose an \emph{admissible direction} $e^{\iu\phi}$ for $\sfC$, that is, an element $e^{\iu\phi} \in S^1$ satisfying $e^{\iu\phi}\notin \R_{>0} (u-u')$ for any $u,u' \in \sfC$.  By the Hukuhara-Turrittin theorem (see \cite[Lemma 8.3]{Hertling-Sevenheck:nilpotent}), the above formal decomposition \eqref{eq:formal_decomp} lifts uniquely to an analytic decomposition 
\begin{equation*} 
\Phi_I \colon \QDM(X)_\tau \Bigr |_{I} \cong \bigoplus_{u\in \sfC} 
 e^{u/z} \otimes \cR_u \Bigr|_{I} 
\end{equation*}
over a sector of the form $I = \{z\in \C^\times : |\arg z - \phi|<\frac{\pi}{2} + \epsilon\}$ for some $\epsilon>0$. Here we mean by ``lifts''  that the map $\Phi_I$ admits, when expressed in terms of local holomorphic frames of $\QDM(X)_\tau$ and $\cR_u$ around $z=0$, an asymptotic expansion as $z\to 0$ along the sector $I$ and that the expansion coincides with $\hPhi$. 

Let $\cS_I$ denote the space of $\nabla$-flat sections over the angular sector $\{\tau\}\times I$: it can be identified with $\cS(X)$ from \S\ref{subsec:Gamma_class} once we specify a lift of the sector $I$ to the universal cover of $\C^\times$.  
The analytic decomposition $\Phi_I$ induces a decomposition of $\cS_I$  
\begin{equation}
\label{eq:decomp_flat_sections} 
\cS_I = \bigoplus_{u\in \sfC} \cV_u  
\end{equation} 
where $\cV_u$ can be identified with the space of flat sections of $\cR_u$ over $I$. 
Since the analytic decomposition $\Phi_I$ is valid over a sector of angle greater than $\pi$, it follows easily that the decomposition \eqref{eq:decomp_flat_sections} is semiorthogonal in the sense that 
\[
[\cV_u, \cV_{u'}) = 0 \quad 
\text{if $\Im(e^{-\iu\phi} u)<\Im(e^{-\iu\phi} u')$}
\]
where $[\cdot,\cdot)$ is the pairing on $\cS_I \cong \cS(X)$ introduced in \eqref{eq:[)}. The data of the vector space $\cS_I$ equipped with the pairing $[\cdot,\cdot)$ and the semiorthogonal decomposition (SOD) \eqref{eq:decomp_flat_sections} constitute a \emph{mutation system} in the sense of \cite[Definition 2.30]{Sanda-Shamoto:analogue}. In what follows, we ignore the torsion part of the $K$-group and write $K(X)$ for $K(X)/\textrm{tors}$. 

\begin{conjecture}[\cite{Sanda-Shamoto:analogue,Galkin:private}]
\label{conj:Gamma_III} 
Suppose that the quantum connection is of exponential type at $\tau\in H^*(X)$. With notation as above, the SOD \eqref{eq:decomp_flat_sections} is induced from a decomposition of the topological $K$-group lattice, i.e.~there exists a decomposition 
\begin{equation} 
\label{eq:decomp_K_group} 
K(X) = \bigoplus_{u\in \sfC} V_u^\phi
\end{equation} 
such that $\cV_u = \frs(V_u^\phi)\otimes \C$, where we identify $\cS_I$ with $\cS(X)$ by choosing a lift $\phi\in \R$ of the direction $e^{\iu\phi} \in I$. (A different choice of the lift $\phi$ changes $V_u^\phi$ by monodromy, i.e.~$V_u^{\phi+2\pi} = V_u^\phi \otimes \omega_X[n]$.)  
\end{conjecture} 

When this conjecture holds, the lattices $\{V_u^\phi\}$ are semiorthogonal with respect to the Euler pairing, i.e.~ $\chi(V_u^\phi, V_{u'}^\phi) = 0$ for $\Im(e^{-\iu\phi} u)<\Im(e^{-\iu\phi} u')$ and therefore the Euler pairing on each $V_u^\phi$ is necessarily unimodular (because the Euler pairing on $K(X)$ is unimodular by Poincar\'e duality). In the semisimple case, we must have $V_u^\phi \cong \Z$ and a generator $\cE$ of $V_u^\phi$ must satisfy $\chi(\cE,\cE) = \pm 1$; the $\hGamma$-conjecture II (Conjecture \ref{conj:Gamma_II}) additionally asserts that $\chi(\cE,\cE)=1$ (this point does not follow from Conjecture \ref{conj:Gamma_III}). 

\begin{remark}
The original formulation in \cite{Sanda-Shamoto:analogue} assumes that $X$ is a smooth Fano variety and claims also that the semiorthogonal decomposition (SOD) \eqref{eq:decomp_K_group} arises from an SOD of the derived category of coherent sheaves. We note that an SOD of the derived category induces an SOD of the topological $K$-group, since projections to the SOD summands are given by Fourier-Mukai kernels in $D^b(X\times X)$ and these kernels induce projections in the topological $K$-group (see the discussion in \cite[\S 4]{Gorchinskiy-Orlov} in the context of algebraic $K$-theory). 
\end{remark}

\begin{example}[\cite{Sanda-Shamoto:analogue}] 
\label{exa:Sanda-Shamoto} 
Sanda and Shamoto proved their conjecture for Fano complete intersections in the projective spaces of Fano index greater than 1. Let $X$ be a degree $d$ Fano hypersurface in $\PP^n$, with $n-d>0$. The set of eigenvalues of the quantum multiplication $(E\star_0) = (c_1(X)\star_0)$ is $\{0\} \cup \{T \zeta : \zeta^{n+1-d} = 1\}$ where $T = (n+1-d) \cdot d^{d/(n+1-d)}$. The multiplicity of $T\zeta$ is one and that of $0$ is the dimension of the primitive cohomology plus $d-1$. In this case, the decomposition \eqref{eq:decomp_flat_sections} at $\tau=0$ arises from (up to mutation) the following SOD of the derived category 
\[
D^b(X) =\langle \cA, \cO, \cO(1),\dots,\cO(n-d)\rangle 
\] 
where $\cO,\dots,\cO(n-k)$ are exceptional objects corresponding to simple eigenvalues $T\zeta$ and $\cA$ is the right orthogonal of $\langle \cO,\cO(1),\dots,\cO(n-d)\rangle$ corresponding to the eigenvalue $0$. 
\end{example} 

The following problem naturally arises: 

\begin{problem} 
\label{prob:piece}
Understand a geometric meaning of each regular singular piece $\cR_u$ and the corresponding unimodular lattice $V_u^\phi$ predicted in Conjecture \ref{conj:Gamma_III}. 
\end{problem}  
In the semisimple case, each regular singular piece is the quantum connection of a point and the $K$-class $\cE_i^\phi$ in the $\hGamma$-conjecture II (Conjecture \ref{conj:Gamma_II})  corresponds to a generator of $K^0(\pt)\cong \Z$. The subcategory $\cA$ in Example \ref{exa:Sanda-Shamoto} is equivalent to the category of graded matrix factorizations of a degree $d$ polynomial $F(x_0,\dots,x_n)$ defining the hypersurface \cite{Orlov:sing}: it is known to be a fractional Calabi-Yau category (in the sense that a  power of the Serre functor equals the shift functor).

\subsection{Monodromy data and Riemann-Hilbert problem} 
\label{subsec:RH} 
Let us assume that $X$ satisfies Conjecture \ref{conj:Gamma_III}. In this section we explain how the SOD \eqref{eq:decomp_K_group} encodes the irregular monodromy (Stokes) data, following \cite[\S 8]{Hertling-Sevenheck:nilpotent} and \cite{Sanda-Shamoto:analogue}. We also formulate a  Riemann-Hilbert problem that reconstructs quantum cohomology from the SOD \eqref{eq:decomp_K_group}, formal data \eqref{eq:formal_decomp} and certain additional data.  

\paragraph{\textbf{Monodromy}} 
The monodromy transformation $T\colon s(z) \mapsto s(e^{2\pi\iu} z)$ on $\cS_I$ is determined from the pairing $[\cdot,\cdot)$ as
\[
[T s_1, s_2) =(s_1(e^{\pi\iu}z), s_2(z)) = (s_2(e^{-\pi\iu}z), s_1(z)) =  [s_2,s_1).  
\]
The restriction $[\cdot,\cdot)_u$ of the pairing $[\cdot,\cdot)$ to $\cV_u$ is non-degenerate and is induced from the pairing $(\cdot,\cdot)_u$ on $\cR_u$. The monodromy transformation $T_u\colon \cV_u \to \cV_u$ on flat sections of $\cR_u$ is likewise determined by $[T_u s_1, s_2)_u = [s_2,s_1)_u$.

\begin{figure}[t]
\centering 
\includegraphics[width=3.7cm]{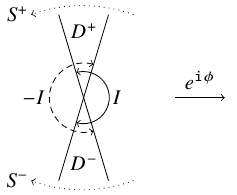}
\caption{The angular sectors $\pm I$ and the paths of analytic continuations used to define the Stokes maps $S^\pm$}
\label{fig:Stokes} 
\end{figure}

\paragraph{\textbf{Stokes data}} Let $\cS_{-I}$ denote the space of $\nabla$-flat sections over the opposite sector $\{\tau\}\times (-I)$. The Poincar\'e pairing $(\cdot,\cdot) \colon \cS_{-I}\times \cS_I\to \C$ identifies $\cS_{-I}$ with the dual space of $\cS_I$ and the decomposition $\cS_{-I} = \bigoplus_{u\in \sfC} \cV'_u$ associated with the sector $-I$ is dual to that for $I$, i.e.~$(\cV'_{u'}, \cV_u) = 0$ for $u\neq u'$. The Stokes data are given by the analytic continuation maps $S^{\pm} \colon \cS_I \to \cS_{-I}$, $s(z) \mapsto s(e^{\pm \pi\iu}z)$. By the very definition of the pairing $[\cdot,\cdot)$, they are determined from $[\cdot,\cdot)$ as 
\[
(S^+s_1,s_2) = [s_2,s_1), \qquad (S^- s_1,s_2) = [s_1,s_2).  
\]
Then we have $T= (S^-)^{-1} S^+$. 
The Stokes maps $S^\pm$ are upper (or lower) triangular in the sense that $S^+(\cV_u) \subset \bigoplus_{\Im(e^{-\iu\phi}u')\ge 
\Im(e^{-\iu\phi}u)} \cV'_{u'}$ 
and $S^-(\cV_u) \subset \bigoplus_{\Im(e^{-\iu\phi}u')\le \Im(e^{-\iu\phi}u)} \cV'_{u'}$. They can be used to glue the connections over the opposite sectors 
\[
\bigoplus_{u\in \sfC}  e^{u/z} \otimes \cR_u\Bigr|_{-I} \quad \text{and} \quad 
\bigoplus_{u\in \sfC} e^{u/z} \otimes \cR_u\Bigr|_{I}
\]
along the two overlapping domains $D^\pm=I \cap (-I) \cap \{\pm \Im(z e^{-\iu\phi})>0\}$ (see Figure \ref{fig:Stokes}). Hence the Stokes data reconstruct an analytic germ of the quantum connection at $z=0$ from the formal data $\{\cR_u\}_{u\in \sfC}$. 

\paragraph{\textbf{Riemann-Hilbert problem}} The global quantum connection over $\PP^1$ can be reconstructed by gluing the germ of the connection at $z=0$ with a connection around $z=\infty$ via the $\hGamma$-integral structure. The quantum connection around $z=\infty$ is gauge-equivalent, via $L(\tau,z)$, to the connection
\[
\nabla^{(\infty)}_{z\partial_z} = z\parfrac{}{z}  - \frac{c_1(X)}{z} + \mu 
\]
on the trivial bundle $F_\infty = H^*(X)\times (\PP^1\setminus \{0\}) \to \PP^1\setminus \{0\}$. We identify the space of $\nabla^{(\infty)}$-flat sections with the $K$-group via the framing $\Psi_\infty \colon K(X) \to H^*(X) \otimes \cO_{\widetilde{\C^\times}}$ (cf.~\eqref{eq:framing}) given by  
\begin{equation} 
\label{eq:Psi_infty}
\Psi_\infty(\alpha): = (2\pi)^{-n/2} z^{-\mu} z^{c_1(X)} \hGamma_X (2\pi\iu)^{\deg/2} \ch(\alpha).
\end{equation} 
We glue the bundle $(F_\infty,\nabla^{(\infty)})$ with the germ around $z=0$ by identifying the flat section $\Psi_\infty(\alpha)$ with $\alpha \in V_u^\phi$ with the flat section in $\cV_u \cong \Gamma(I,\cR_u)^\nabla$ corresponding to $\alpha$ (here we need an identification $V_u^\phi \cong \cV_u$). This gives us a global vector bundle $\hF \to \PP^1$ with a meromorphic connection $\hnabla$. The glued bundle $\hF$ must be trivial (although it is not a priori clear); the trivialization of $F_\infty$ at $z=\infty$ induces a trivialization $\hF \cong H^*(X) \times \PP^1$. The pair $(\hF,\hnabla)$ is identified with the quantum connection at $\tau$. 

More explicitly, this reconstruction procedure can be described as the following Riemann-Hilbert problem for functions $Y_\pm=(\Phi_{\pm I})^{-1}$ (over the sectors $\pm I$) and $Y_\infty = L(\tau,z)$ (around $z=\infty$). This is an extension of the Riemann-Hilbert problem described by Dubrovin \cite{Dubrovin:2D}, \cite[Lecture 4]{Dubrovin:Painleve} in the semisimple case. 

\begin{problem} 
\label{prob:RH} 
Suppose that we are given the following data $\sfC$, $\{\cR_u\}_{u\in \sfC}$, $e^{\iu\phi}$, $I$, $\phi$, $(K(X),\chi)$,  $K(X)= \bigoplus_{u\in \sfC} V_u$, $\Psi_u$, $\Psi_\infty$: 
\begin{itemize} 
\item a subset $\sfC$ of $\C$;   
\item a finite free $\C\{z\}$-module $\cR_u$ with a regular singular connection for each $u\in \sfC$; 
\item an admissible direction $e^{\iu\phi}$ for $\sfC$ and a sector $I = \{z\in \C^\times : |\arg z - \phi|<\frac{\pi}{2}+\epsilon\}$ centred around it;   
\item  a unimodular lattice $(K(X),\chi)$ of rank $\dim H^*(X)$, a lift $\phi\in \R$ of $e^{\iu\phi}$ and an SOD $K(X)= \bigoplus_{u\in \sfC} V_u$; 
\item a framing $\Psi_u \colon V_u \to \Gamma(I,\cR_u)^\nabla$ for each $u\in \sfC$ such that $\Psi_u$ induces an isomorphism over $\C$ and intertwines the transformation $T_u \in \End(V_u)$ given by $\chi(T_u \alpha,\beta) = \chi(\beta,\alpha)$ with the monodromy $s(z) \mapsto s(e^{2\pi \iu} z)$ on $\Gamma(I,\cR_u)^\nabla$;  
\item the ``$\hGamma$-integral'' framing $\Psi_\infty\colon K(X) \to H^*(X)\otimes \cO_{\widetilde{\C^\times}}$ given in \eqref{eq:Psi_infty}, which satisfies $\Psi_\infty(e^{2\pi \iu}z) = \Psi_\infty(z)\circ T$ with $T\in \End(K(X))$ given by $\chi(T \alpha,\beta) = \chi(\beta,\alpha)$. 
\end{itemize} 
We define Stokes maps $S^\pm \colon K(X) \to K(X)^\vee$ by $\langle S^+ \alpha, \beta\rangle =\chi(\beta,\alpha)$, $\langle S^-\alpha,\beta \rangle = \chi(\alpha,\beta)$ and a framing $\Psi_{-,u}\colon V_u^\vee \to \Gamma(-I,\cR_u)^\nabla$ over the opposite sector $-I$ by 
\[
\Psi_{-,u}(\chi(\alpha,\cdot)) := \text{clockwise analytic continuation of $\Psi_u(\alpha)$ through $D_-$} 
\]
for $\alpha\in V_u$. 
We set 
\begin{align*} 
\Psi & := \oplus_{u\in \sfC} \Psi_u \colon K(X) =\oplus_{u\in \sfC} V_u \to \oplus_{u\in \sfC} \Gamma(I,\cR_u)^\nabla \\ 
\Psi_- & := \oplus_{u\in \sfC} \Psi_{-,u}\colon K(X)^\vee = \oplus_{u\in \sfC} V_u^\vee \to \oplus_{u\in \sfC} \Gamma(-I,\cR_u)^\nabla
\end{align*} 
The problem is to find (matrix-valued) holomorphic functions 
\[
Y_\infty \in GL(H^*(X))\otimes \cO_{\PP^1\setminus \{0\}}, \quad 
Y_\pm \colon \oplus_{u\in \sfC} \cR_u\bigr|_{\pm I} \to H^*(X) \otimes \cO_{\pm I} 
\]
such that 
\begin{align*} 
Y_\infty|_{z=\infty} = \id, \qquad Y_\pm  \to Y_0 \quad \text{as $z\to 0$ along the sector $\pm I$} 
\end{align*} 
for an invertible operator $Y_0\colon \oplus_{u\in \sfC} \cR_u|_{z=0} \to H^*(X)$ and that  
\begin{align*} 
Y_+ \Psi e^{-U/z} & = Y_\infty \Psi_\infty  && \text{over $I$} \\
Y_- \Psi_- e^{-U^\vee/z} S^\pm & = Y_+ \Psi e^{-U/z} && \text{over $D^\pm$}   
\end{align*} 
where $D^\pm$ is as before, the determination of $\Psi_\infty$ over $I$ is given by $|\arg z - \phi|<\frac{\pi}{2}+\epsilon$, $U:= \oplus_{u\in 
\sfC} u \id_{V_u} \in \End(K(X))$ and $U^\vee := \oplus_{u\in \sfC} u \id_{V_u^\vee} \in \End(K(X)^\vee)$. 
\end{problem} 

A solution $(Y_\pm, Y_\infty)$ to this problem is unique if exists. The solution $Y_\infty$ gives the fundamental solution $L(\tau,z)$ and hence recovers the quantum connection. It is interesting to note that we reconstruct not only the connection but also the fundamental solution $L(\tau,z)$ (called \emph{calibration} in the theory of Frobenius manifolds): this implies that the value of the parameter $\tau$ can be reconstructed by the asymptotics $L(\tau,z)^{-1} 1 = 1 + \tau z^{-1} + O(z^{-2})$ if we know the unit class $1$. 

\begin{remark} 
The additional data we need here (other than those we already mentioned) is the framing $\Psi_u$ for each regular singular piece. In the semisimple case, we have $\cR_u\cong (\C\{z\},d)$ and $V_u \cong \Z$, so there is essentially a unique choice for $\Psi_u$. A natural candidate for $\Psi_u$ could be given by answering Problem \ref{prob:piece}. See \S \ref{subsec:blowup} for the example where we have a natural candidate for the framing. 
\end{remark} 

\begin{remark} 
If we include odd classes, the monodromy transformation $T$ on $\cS_I$ is given by $(-1)^{\deg \alpha} [T\alpha,\beta) = [\beta,\alpha)$; the Stokes maps $S^\pm \colon \cS_I \to \cS_{-I}$ are given by $(S^+ \alpha,\beta) = (-1)^{\deg \alpha} [\beta,\alpha)$, $(S^-\alpha,\beta) = [\alpha,\beta)$. Problem \ref{prob:RH} can be also modified accordingly, using the fact that the pairing $[\cdot,\cdot)$ on $\cS_I$ corresponds to $-\iu \chi(\cdot,\cdot)$ on $K^1(X)$ (see Remark \ref{rem:odd}). 
\end{remark}

\section{Functoriality of quantum cohomology} 
In this section, we discuss a conjectural functoriality of quantum cohomology under birational transformations. Roughly speaking, we expect that the relationship between quantum cohomology is induced from a natural map between $K$-groups via the $\hGamma$-integral structure. Let $X_1$, $X_2$ be smooth projective varieties and let $\varphi \colon X_1 \dasharrow X_2$ be a birational map. Suppose that $\varphi$ fits into the following commutative diagram 
\begin{equation} 
\label{eq:birational} 
\xymatrix{
& \hX \ar[ld]_{p_1}\ar[rd]^{p_2} & \\
X_1 \ar@{-->}[rr]^{\varphi}& & X_2.
}
\end{equation} 
where $p_1, p_2$ are projective birational morphisms. We say that $\varphi$ is \emph{crepant} (or $K$-equivalent) if $p_1^*K_{X_1} = p_2^*K_{X_2}$ and \emph{discrepant} otherwise. We allow $X_i$ to be smooth Deligne-Mumford stacks (with projective coarse moduli spaces) so that we can include crepant resolutions of orbifolds in the following discussion. 

\subsection{Crepant transformation} 
Suppose that $\varphi\colon X_1 \dasharrow X_2$ is crepant. In this case it can be shown that $H^*(X_1) \cong H^*(X_2)$ as graded vector spaces by Kontsevich's motivic integration (see e.g.~\cite{Yasuda}). A famous conjecture of Yongbin Ruan \cite{Ruan:crepant} says that the quantum cohomologies of $X_1$ and $X_2$ become isomorphic after analytic continuation. This problem has been studied by many people, see e.g.~\cite{Li-Ruan, Bryan-Graber, LLW:flop,McLean:birat}. We give a version of the conjecture stated in terms of quantum D-modules and the $\hGamma$-integral structure following \cite[Conjecture 5.1]{CIT}, \cite[\S 5.5]{Iritani:Int}, \cite{Iritani:Ruan, CR:conjecture}. 
\begin{conjecture}[Crepant Transformation Conjecture] 
\label{conj:CTC} 
Let $\varphi \colon X_1 \dasharrow X_2$ be a crepant birational map. 
There exists a map $f$ from an open subset of $H^*(X_1)$ to an open subset of $H^*(X_2)$ such that, after analytic continuation, we have an isomorphism of quantum D-modules $\QDM(X_1) \cong f^* \QDM(X_2)$. Moreover,via the $\hGamma$-integral structure, the isomorphism is induced by an isomorphism $(K(X_1),\chi) \cong (K(X_2),\chi)$ of topological $K$-group lattices. 
\end{conjecture} 

Recall from \S\ref{subsec:QDM} that the quantum D-module $\QDM(X_i)$ is the tuple of the cohomology bundle $F$, the quantum connection and the Poincar\'e pairing and the isomorphism in the conjecture is required to respect these structures. Conjecture \ref{conj:CTC} was proved\footnote{For complete intersections, we restrict to the ambient part of quantum cohomology in \cite{CIJ}.} for crepant transformations between complete intersections in toric Deligne-Mumford stacks, which are induced from variation of GIT quotients \cite{CIJ}. In that case, it is shown that the map $K(X_1) \cong K(X_2)$ between $K$-groups is given by a Fourier-Mukai transformation that gives rise to the equivalence of derived categories of $X_1$ and $X_2$. The calculation needed in this result is an extension of the work of Borisov-Horja \cite{Borisov-Horja:FM} that relates analytic continuation of hypergeometric solutions to the GKZ system and Fourier-Mukai transformations between toric orbifolds.  
\begin{remark} 
(1) We can hope that the isomorphism $K(X_1) \cong K(X_2)$ is induced by an equivalence of derived categories in general. A different derived equivalence can arise from a different choice of paths of analytic continuation. 

(2) When Conjecture \ref{conj:CTC} holds, the map $f$ is necessarily a local isomorphism and identifies the $F$-manifold structure \cite{HM} of quantum cohomology. In the case of crepant resolutions of orbifolds, it has been observed in \cite{CIT} that $f$ is not necessarily affine-linear unless the orbifold satisfies the Hard Lefschetz condition. 
\end{remark} 

\subsection{Discrepant transformation} 
We present a conjectural picture in the discrepant case following \cite{Iritani:discrepant}. In the discrepant case, the ranks of cohomology are different in general and we expect to have an \emph{orthogonal} decomposition of formal quantum D-modules and a \emph{semiorthogonal} decomposition of the $\hGamma$-integral structure. As in \S\ref{subsec:Gamma_III}, we write $\ovQDM(X) := \QDM(X) \otimes_{\cO[z]} \cO[\![z]\!]$ for the quantum D-module formalized along $z=0$. Because of the lack of abundant evidences, we state our picture as problems rather than conjectures. 
\begin{problem}[formal decomposition] 
\label{prob:discrepant_formal} 
Let $\varphi \colon X_1 \dasharrow X_2$ be a birational map fitting into the diagram  \eqref{eq:birational} such that $p_1^*K_{X_1} - p_2^*K_{X_2}$ is an effective divisor. Show that there exists a map $f$ from an open subset of $H^*(X_1)$ to an open subset of $H^*(X_2)$ such that we have an orthogonal decomposition 
\[
\ovQDM(X_1) \cong f^* \ovQDM(X_2) \oplus \scrD
\] 
where $\scrD$ is a locally free $\cO[\![z]\!]$-module equipped with a flat meromorphic connection $\nabla^{\scrD}$ and a $\nabla^{\scrD}$-flat pairing $(\cdot,\cdot)_{\scrD}\colon (-)^* \scrD \otimes \scrD \to \cO[\![z]\!]$. 
\end{problem} 

Problem \ref{prob:discrepant_formal} has been solved for discrepant birational transformations between toric Deligne-Mumford stacks which arise from a variation of GIT quotients \cite{Iritani:discrepant}. The proof is based on mirror symmetry for toric stacks \cite{CCIT:mirrorthm, CCIT:MS}. 

Suppose that Problem \ref{prob:discrepant_formal} is solved for some $\varphi \colon X_1 \dasharrow X_2$, and also suppose (for simplicity) that $\ovQDM(X_1)$ is of exponential type (see \S\ref{subsec:Gamma_III}) at some $\tau\in H^*(X_1)$ in the domain of the map $f$. Then $\ovQDM(X_2)_{f(\tau)}$ and $\scrD|_\tau$ are also of exponential type. We further assume the following: there exist a phase $\phi \in \R$ and real numbers $l_1>l_2$ such that 
\begin{itemize} 
\item every eigenvalue $u$ of $-\nabla^{\scrD}_{z^2 \partial_z} \in \End(\scrD|_{z=0,\tau})$ satisfies either $\Im(e^{-\iu\phi}u)>l_1$ or $l_2>\Im(e^{-\iu\phi}u)$ and 
\item every eigenvalue $u$ of $(E^{X_2}\star_{f(\tau)})\in \End(H^*(X_2))$ satisfies $l_1>\Im(e^{-\iu\phi}u)>l_2$. 
\end{itemize} 
Then $\scrD|_\tau$ decomposes as $\scrD|_\tau  = \scrD_1 \oplus \scrD_2$ so that every eigenvalue of $-\nabla^{\scrD_1}_{z^2 \partial_z}$ on $\scrD_1|_{z=0}$ satisfies $\Im(e^{-\iu\phi} u) >l_1$ and that every eigenvalue of $-\nabla^{\scrD_2}_{z^2 \partial_z}$ on $\scrD_2|_{z=0}$ satisfies $\Im(e^{-\iu\phi} u) <l_2$. By varying $\phi$ a little, we may assume that $e^{\iu\phi}$ is admissible for the eigenvalues of $(E^{X_1}\star_\tau)$.  As discussed in \S\ref{subsec:Gamma_III}, by the Hukuhara-Turrittin theorem, the formal decomposition $\ovQDM(X_1)_\tau \cong \scrD_1\oplus \ovQDM(X_2)_{f(\tau)} \oplus \scrD_2$ lifts to an analytic decomposition of connections over a sector of the form $I=\{z \in \C^\times: |\arg z - \phi|<\frac{\pi}{2} + \epsilon \}$ for some $\epsilon>0$ 
\begin{equation} 
\label{eq:analytic_decomp} 
\QDM(X_1)_\tau\Bigr|_I \cong \scrD_{1,I} \oplus \QDM(X_2)_\tau\Bigr|_I \oplus \scrD_{2,I}
\end{equation} 
where $\scrD_{i,I}$ is an analytic connection over the sector $I$. 
\begin{problem}[analytic decomposition]  
\label{prob:discrepant_analytic} 
Show that the analytic decomposition \eqref{eq:analytic_decomp} is induced, via the $\hGamma$-integral structures for $X_1$ and $X_2$, by an SOD of topological $K$-groups: 
\begin{equation} 
\label{eq:SOD_K-group} 
K(X_1) \cong K_1 \oplus K(X_2) \oplus K_2 
\end{equation} 
such that the associated inclusion $K(X_2) \to K(X_1)$ respects the Euler pairing. 
\end{problem} 

Problem \ref{prob:discrepant_analytic} has been answered affirmatively when $X_1$, $X_2$ are weak-Fano compact toric Deligne-Mumford stacks (satisfying certain mild technical conditions) and $\varphi \colon X_1 \to X_2$ is a weighted blowup (or a root construction) along a toric substack $Z$ \cite{Iritani:discrepant}. We also showed that the decomposition \eqref{eq:SOD_K-group} at some $\tau$ and $\phi$ is given by an Orlov-type SOD \cite{Orlov:projective}. 
We could hope that the SOD \eqref{eq:SOD_K-group} in $K$-theory arises from an SOD $D^b(X_1) \cong \langle \cA_1, D^b(X_2), \cA_2 \rangle$ of the derived category; such an SOD in the derived category has been conjectured in \cite{BFK:VGIT}. 
\begin{remark}
There are closely related works by Bayer \cite{Bayer:blowup}, Acosta-Shoemaker \cite{Acosta-Shoemaker:toric,Acosta-Shoemaker:blowup_LG} and Gonz\'alez-Woodward \cite{GW:tmmp}. 
A formal decomposition of quantum D-modules under flips similar to our picture has been also proposed by Lee, Lin and Wang \cite{LLW:proc, LLW:flips}. 
\end{remark} 

\begin{remark} 
In Problem \ref{prob:discrepant_formal}, $f$ is necessarily a submersion and the  quantum cohomology $F$-manifold of $X_1$ locally decomposes into the product of the quantum cohomology $F$-manifold of $X_2$ and that corresponding to $\scrD$. \emph{Proof}. The map $T_\tau H^*(X_1) \to \ovQDM(X_1)|_{z=0,\tau} = H^*(X_1)$ given by $v\mapsto z \nabla_v 1$ is an isomorphism. If the unit section $1$ maps to $(s,t) \in \ovQDM(X_2)\oplus \scrD$ under the isomorphism, it follows that the map $df(T_\tau H^*(X_1)) \to \ovQDM(X_2)|_{z=0,f(\tau)}$, $w \mapsto z \nabla_w s$ is surjective. This can only happen when $df\colon T_\tau H^*(X_1) \to T_{f(\tau)} H^*(X_2)$ is surjective. The ring homomorphism $T_\tau H^*(X_1) \hookrightarrow \End(\ovQDM(X_1)|_{z=0,\tau})$, $v \mapsto z \nabla_v$ factors through $T_\tau H^*(X_1) \to T_{f(\tau)} H^*(X_2) \oplus \End(\scrD|_{z=0,\tau})$ and this gives a decomposition of the ring $(T_\tau H^*(X_1),\star_\tau)$. The flatness of $\nabla$ shows that the decomposition is integrable. 
\end{remark} 

\begin{remark} 
For a higher-genus generalization of Conjecture \ref{conj:CTC} and Problem \ref{prob:discrepant_formal}, we refer the reader to \cite{CI:Fock}, \cite{Iritani:discrepant}. 
\end{remark} 

\subsection{Riemann-Hilbert problem for blowups} 
\label{subsec:blowup} 
Let $X$ be a smooth projective variety and let $Z \subset X$ be a smooth subvariety. Let $\varphi \colon \tX \to X$ be the blowup of $X$ along $Z$. The above mentioned results for toric blowups suggest the following conjectural reconstruction algorithm for quantum cohomology of $\tX$ from quantum cohomology of $X$ and $Z$. This is similar to the procedure in \S\ref{subsec:RH}. 

\paragraph{\textbf{Orlov decomposition}} Let $c$ be the codimension of $Z$ in $X$. By Orlov \cite{Orlov:projective} we have the SOD of the $K$-group: 
\begin{align*}
K(\tX) & =  \varphi^*K(X) \oplus K(Z)_0 \oplus \cdots \oplus K(Z)_{c-2} 
\cong K(X) \oplus K(Z)^{\oplus (c-1)} 
\end{align*} 
where $K(Z)_k = j_*(\cO(k) \otimes \pi^* K(Z))$ with $j \colon E\hookrightarrow X$ the inclusion of the exceptional locus and $\pi \colon E \cong \PP(N_{Z/X}) \to Z$ a projective bundle. We shall fix this decomposition. The cohomology of $\tX$ is isomorphic to $H^*(X) \oplus H^{*-2}(Z) \oplus \cdots \oplus H^{*-2c+2}(Z)$ as graded vector spaces. The cup product structure on $H^*(\tX)$, the $\hGamma$-class and the Chern character for $\tX$ can be reconstructed from those for $X$, $Z$, the push-forward and pull-back maps between $H^*(X)$, $H^*(Z)$ and the Chern classes $c_i(N_{Z/X}) \in H^{2i}(Z)$. 

\paragraph{\textbf{Formal data}} We choose parameters $\sigma\in H^*(X)$ and $\rho_0, \dots, \rho_{c-2} \in H^*(Z)$ and a phase $\phi \in \R$ so that $\Im(e^{-\iu\phi} v)> \Im(e^{-\iu\phi} u_0) > \cdots >\Im(e^{-\iu \phi} u_{c-2})$ for all eigenvalues $v$ of $(E^X \star_\sigma)$ and all eigenvalues $u_i$ of $(E^Z\star_{\rho_i})$. We define $\ovQDM:=\ovQDM(X)_\sigma \oplus \ovQDM(Z)_{\rho_0} \oplus \cdots \oplus \ovQDM(Z)_{\rho_{c-2}}$. This will be the formal quantum D-module for $\tX$. 

\paragraph{\textbf{Gluing}} 
The given formal decomposition for $\ovQDM$ should lift to analytic decompositions over the sectors $I$ and $-I$, with $I = \{z\in \C^\times : |\arg z - \phi|<\frac{\pi}{2}+\epsilon\}$ for some $\epsilon>0$ 
\[
\Phi_{\pm I} \colon \QDM\bigr|_{\pm I} \cong \QDM(X)_\sigma \oplus \QDM(Z)_{\rho_0} \oplus \cdots \oplus \QDM(Z)_{\rho_{c-2}} \bigr|_{\pm I} 
\]
and the two analytic decompositions should be glued together by the Stokes data induced from Orlov's SOD. Finally we glue it with the connection near $z=\infty$ via the $\hGamma$-integral structure to get the quantum D-module for $\tX$. 

The reconstruction can be formulated as a Riemann-Hilbert problem for $Y_\pm=(\Phi_{\pm I})^{-1}$ and $Y_\infty=L(\tau,z)$ (a fundamental solution for $\tX$, see \eqref{eq:fundsol}) as follows. Define $S^\pm \colon K(\tX) \to K(\tX)^\vee$ by $\langle S^+ \alpha,\beta\rangle = \chi(\beta,\alpha)$, $\langle S^- \alpha,\beta\rangle = \chi(\alpha,\beta)$ as before; also define $\Psi\colon K(\tX) \cong K(X)\oplus K(Z)^{\oplus (c-1)} \to (H^*(X)\oplus H^*(Z)^{\oplus (c-1)})\otimes \cO_I$ as 
\[
\Psi(\alpha,\beta_0,\dots,\beta_{c-2}) =\frs_{X}(\alpha) (\sigma,z)\oplus \frs_Z(\beta_0)(\rho_0,z) \oplus \cdots \oplus \frs_Z(\beta_{c-2})(\rho_{c-2},z) 
\] 
where $\frs_X$, $\frs_Z$ are the maps \eqref{eq:framing} defined for $X$ and $Z$ respectively and define $\Psi_- \colon K(\tX)^\vee \cong K(X)^\vee \oplus (K(Z)^\vee)^{\oplus (c-1)} \to (H^*(X)\oplus H^*(Z)^{\oplus (c-1)})\otimes \cO_{-I}$ as
\begin{align*}
\Psi_-(\chi(\alpha,\cdot),\chi(\beta_0,\cdot),\dots,\chi(\beta_{c-2},\cdot)) 
= \text{clockwise analytic continuation of $\Psi(\alpha,\beta_0,\dots,\beta_{c-2})$.}
\end{align*} 
Let $\Psi_\infty$ be the map \eqref{eq:Psi_infty} with $X$ there replaced with $\tX$. The problem is to find functions 
\[
Y_\infty \in GL(H^*(\tX)) \otimes \cO_{\PP^1\setminus \{0\}}, \quad Y_\pm \in \Hom(H^*(X)\oplus H^*(Z)^{\oplus (c-1)},H^*(\tX))\otimes \cO_{\pm I}. 
\] 
such that $Y_\infty|_{z=\infty} = \id$, $Y_\pm \to Y_0$ as $z\to 0$ along the sector $\pm I$ for an invertible operator $Y_0$ and that 
\begin{align*} 
Y_+ \Psi  & = Y_\infty \Psi_\infty && \text{over $I$}\\
Y_- \Psi_- S^\pm & = Y_+ \Psi && \text{over $D^\pm$} 
\end{align*} 
where $D^\pm$ is as before. As discussed in \S\ref{subsec:RH}, we can reconstruct the value of the parameter $\tau$ for the quantum D-module of $\tX$ and it becomes a function of $\sigma,\rho_0,\dots,\rho_{c-2}$; the parameter space locally splits into the product of $H^*(X)$ and $(c-1)$ copies of $H^*(Z)$ as an $F$-manifold. 

\begin{remark} 
Recently, Katzarkov, Kontsevich and Pantev \cite{Kontsevich:talk} formulated a closely related conjecture for quantum cohomology of blowups and gave a remarkable application to the problem of rationality.   
\end{remark} 


\paragraph{\textbf{Acknowledgements}}
I thank Sergey Galkin for valuable comments on a draft version of the paper and allowing me to present his formulation in this paper. 

\paragraph{\textbf{Funding}} 
This work was partially supported by JSPS grant 16H06335, 20K03582 and 21H04994.

\end{document}